\documentclass[11pt,letterpaper]{amsart}
\usepackage{amsmath,amsfonts,amssymb,amsthm}
\usepackage{hyperref}
\usepackage{mathtools}
 \usepackage{pdfpages}
 \usepackage{cite}

\newtheorem{lem}{Lemma}[section]
\newtheorem{teo}[lem]{Theorem}
\newtheorem{thm}[lem]{Theorem}
\newtheorem{pro}[lem]{Proposition}
\newtheorem{cor}[lem]{Corollary}
\newtheorem{claim}[lem]{Claim}
\newtheorem{defi}[lem]{Definition}

\newtheorem{conj}{Conjecture}

\newtheorem*{rem*}{Remark}
\newtheorem*{teo*}{Theorem}

\newcounter{claimcounter}
\numberwithin{claimcounter}{lem}

\newtheorem{thmA}{Theorem}
\theoremstyle{definition}

\usepackage{scalerel,stackengine}
\stackMath
\newcommand\reallywidehat[1]{%
\savestack{\tmpbox}{\stretchto{%
  \scaleto{%
    \scalerel*[\widthof{\ensuremath{#1}}]{\kern-.6pt\bigwedge\kern-.6pt}%
    {\rule[-\textheight/2]{1ex}{\textheight}}
  }{\textheight}%
}{0.5ex}}%
\stackon[1pt]{#1}{\tmpbox}%
}

\newcommand{\myle}[1]{\ensuremath{\stackrel{\text{#1}}{\leq}}}

\usepackage{stmaryrd}

\usepackage{commath}

\usepackage[all,cmtip]{xy}

\usepackage{tikz-cd}
 \usetikzlibrary{positioning}

\usepackage{mathtools}

\DeclareMathOperator{\Id}{Id}

\DeclareMathOperator{\D}{\mathcal D}

\DeclareMathOperator{\im}{Im}

\DeclareMathOperator{\Hom}{Hom}
\DeclareMathOperator{\Ext}{Ext}

\DeclareMathOperator{\cd}{cd}

\DeclareMathOperator{\Aut}{Aut}

\DeclareMathOperator{\Tor}{Tor}

\renewcommand{\b}{b^{(2)}}

\newcommand{\ga}{\gamma}

\newcommand{\Ga}{\Gamma}

\newcommand{\p}{\F_p}

\newcommand{\hp}{\widehat{\mathcal N_p}}

\newcommand{\FP}{\mathtt{FP}}

\newcommand{\F}{\mathbb{F}}
\newcommand{\PD}{\mathtt{PD}}

\newcommand{\Gc}{G_{\widehat {\mathcal C}}}

\newcommand{\Hcc}{H_{\widehat {\mathcal N_p}}}

\newcommand{\Gsr}{\p \llbracket G_{\widehat {\mathcal S}}\rrbracket}

\newcommand{\Ger}{\p \llbracket G_{\widehat {\mathcal E}}\rrbracket}
\newcommand{\Gcr}{\p \llbracket G_{\widehat {\mathcal C}}\rrbracket}

\newcommand{\Dp}{\D_{\p[G]}}

\newcommand{\Sg}{\Sigma}

\newcommand{\sagecode}[1]{}

\newcommand{\GG}{{\bf G}}

\newcommand{\UU}{{\bf U}}
\newcommand{\VV}{{\bf V}}
\newcommand{\MM}{{\bf M}}

\renewcommand{\le}{	<}
\renewcommand{\geq}{\geqslant}
\renewcommand{\leq}{\leqslant}

\newcommand{\po}{\colon}
\newcommand{\Hc}{ H_{\textrm{cts}}}
\newcommand{\Hcp}{H_{\textrm{cts}}}
\renewcommand{\hat}{\widehat}

\newcommand{\Cor}{\textrm{Cor}}
\newcommand{\Res}{\textrm{Res}}

\newcommand{\Z}{\mathbb{Z}}

\newcommand{\N}{\mathbb{N}}

\newcommand{\Q}{\mathbb{Q}}

\newcommand{\lrar}{\longrightarrow}

\newcommand{\rar}{\rightarrow}

\newcounter{ismaelcomments}

\newcounter{andreicomments}

 \date{\today}

\title[Prosolvable rigidity of surface groups ]{Prosolvable rigidity of surface groups}
\author{Andrei Jaikin-Zapirain}
 
\address{Departamento de Matem\'aticas, Universidad Aut\'onoma de Madrid \and  Instituto de Ciencias Matem\'aticas, CSIC-UAM-UC3M-UCM}
\email{andrei.jaikin@uam.es}

\author{Ismael Morales}
\address{Mathematical Institute, University of Oxford, Radcliffe Observatory, Andrew Wiles Building, Woodstock Rd, Oxford OX2 6GG}
\email{morales@maths.ox.ac.uk}

\begin{document}
\begin{abstract} Surface groups are known to be the Poincar\'e Duality groups of dimension two since the work of Eckmann, Linnell and M\"uller. We prove a prosolvable analogue of this result that allows us to show that surface groups are profinitely (and prosolvably) rigid among finitely generated   groups that satisfy $\cd(G)=2$ and $b_2^{(2)}(G)=0$. We explore two other consequences.

On the one hand, we derive that if $u$ is a surface word  of a finitely generated free group $F$ and $v\in F$ is measure equivalent to $u$ in all finite solvable quotients of $F$ then $u$ and $v$ belong to the same $\Aut(F)$-orbit. 

Finally, we get a partial result towards Mel'nikov's surface group conjecture.   Let $F$ be a free group of rank $n\geq 3$ and let $w\in F$. Suppose that $G=F/\langle\!\langle w\rangle\!\rangle$ is a residually finite  group all of whose finite-index subgroups are one-relator groups. Then  $G$ is 2-free. Moreover, we show that if $H^2(G; \Z)\neq 0$ then $G$ must be a surface group.   
\end{abstract}

 \maketitle

\section{Introduction} \subsection{Profinite rigidity.} A compelling question in group theory that has recently received much attention is  to what extend  a group is determined by its finite quotients. This subject is  known by the name of profinite rigidity and has been explored by a rich variety of techniques ranging from algebra to geometry and topology. One of the most intriguing open questions is whether  finitely generated free groups and surface groups are profinitely rigid. In this article we address this and related questions for  surface groups (being Theorems \ref{main} and \ref{Mel'Intro} our main results). By {\bf surface group} we will always mean the fundamental group of a closed surface.

The criterion we use to recognise surface groups is the well-known fact that these are the Poincar\'e duality groups of cohomological dimension two by the work of Eckmann--Linnell--M\"{u}ller \cite{Eck80, Eck83} (subsequently generalised by Bowditch \cite{Bow04}).   Theorem \ref{main0}  is a prosolvable extension of this principle that will allow us to recognise surface groups from their finite quotients in various situations.  
\begin{thmA}\label{main0} Let $p$ be a prime and let $G$ be a   RFRS group of cohomological dimension two and   type $\FP_2(\F_p)$. Suppose that $G$ is prosolvable $p$-good and that its prosolvable completion $G_{\hat{\mathcal S}}$ is   Poincar\'e Duality of dimension two  at $p$. Then $G$ is a surface group. 
\end{thmA}
The notion of prosolvable $p$-goodness is an analogue of Serre's notion of cohomological goodness \cite{Ser97} that we review in Section \ref{goodsection}. We also recall the notion of Poincar\'e duality in Definition \ref{PD2profinite}. A simple criterion to produce these groups is offered by Kochloukova--Zalesskii \cite{Koc08}. They showed that the pro-solvable completion of an abstract Poincar\'e duality group $\Ga$ of dimension $n$ over $\F_p$ that is $p$-good  is also Poincar\'e duality of dimension $n$ at $p$.  One of the motivations of this article was to get a two-dimensional converse of this fact (as we  do in Theorem \ref{main0}, strengthened in Theorem \ref{mainC0} for other varieties). 

On the other hand,   RFRS groups (which stand for ``residually finite rationally solvable'') form  a class of groups that played a fundamental role in Agol's  solution of the virtual fibring theorem \cite{Ago13} (see also \cite{Kie18}). We refer the reader to Section \ref{finiteSection} for definitions. 
 Theorem \ref{main0} suggests that RFRS groups are also very useful to handle questions on profinite rigidity. In fact,  we can prove the following. 
    
\begin{thmA}\label{main}
Let $S$ be a surface group and let $G$ be a   finitely generated group with $\cd(G)=2$ and $b_2^{(2)}(G)=0$. Assume that either
\begin{enumerate}
\item $G$ is residually finite and has the same finite quotients as $S$; or
\item $G$ is residually-(finite solvable) and has the same finite solvable quotients as $S$.
\end{enumerate}
 Then $G\cong S$.
\end{thmA}

Theorem \ref{main} is applicable to many classes of groups. For example, it includes the class of one-relator groups (because they have $\cd(G)\leq 2$ if they are torsion-free by Lyndon \cite{Lyndon50} and they have $b_2^{(2)}(G)=0$ by Dicks--Linnell \cite{Dic07}). On the other hand, it also includes limit groups for the following reason. If $L$ is a limit group then $b_2^{(2)}(L)=0$ by \cite{BriKoc17}. If, in addition, $L$ has the same prosolvable completion as a non-abelian free or surface group $S$, it follows from \cite[Theorem A]{Mor22}  that $L$ does not contain $\Z^2$ and hence   $\cd(L)\leq 2$ by Sela's finite cyclic hierarchy \cite{Sel01}. So Theorem \ref{main} proves that surface groups can be recognised among limit groups uniquely by their finite solvable quotients, which strengthens previous results obtained in \cite{Wil21, Mor22, Fru23}. The analogue of Theorem \ref{main} for the case when $S$ is free remains open.

\subsection{Measure equivalence.} The fact that Theorem \ref{main} applies to torsion-free one-relator groups allows us to obtain in Theorem \ref{surfaceword} new results in the subject of measure equivalent words. Before discussing them, we recall a few definitions. Let $F$ be a finitely generated free group on free $d$-generators $x_1,\ldots, x_d$ and $w\in F$. Given a finite group $G$, we can define a map $w_G\po G^d\lrar G$ that sends $(g_1,\ldots, g_d)$ to the image of $w$ under the homomorphism $F\lrar G$ sending $x_i$ to $g_i$. We say that two words $w, u\in F$ are {\bf measure equivalent in $G$} if for every  $g\in G$, $|w_G^{-1}(g)|=|u_G^{-1}(g)|$.  If $\mathcal C$ is a non-empty family of finite groups, we say that $w, u\in F$ are {\bf measure equivalent in $\mathcal C$} if they are measure equivalent in every $G\in \mathcal C$. For example if $u=\phi(w)$ for some $\phi\in \Aut(F)$ then $u$ and $w$ are measure equivalent  in the class of all finite groups. The main conjecture in the subject is the following.

\begin{conj}
Let $u,w\in F$ be two words of a finitely generated free group. Assume that $u$ and $w$ are measure equivalent in all finite groups. Then there exists $\phi\in \Aut(F)$ such that $u=\phi(w)$.
\end{conj}
The conjecture is known when one of the words is primitive \cite{PP15} (see also \cite{Wi19,GJ22}) and when one of the words is a surface word \cite{Wi21} (see also \cite{MP21}). We can reprove the conjecture for surface words and, moreover, we show that it is enough to look at the finite solvable quotients. 
\begin{thmA} \label{surfaceword}
Assume that either
\begin{enumerate}
\item $F$ is a free group freely generated by generators $x_1,y_1,\ldots, x_d,y_d$ and  $w=[x_1,y_1]\cdots [x_d,y_d]$ or
\item $F$ is a free group freely generated by generators $x_1,\ldots, x_k$ and $w=x_1^2\cdots x_k^2$.
\end{enumerate}
 Let  $u\in F$. Assume that $w$ and $u$ are measure equivalent in the class of finite solvable groups.  Then there exists $\phi\in \Aut(F)$ such that  $u=\phi(w)$.

\end{thmA}
In fact, we will prove that  it is enough to consider the quotients in the pseudovariety $ \mathcal N_p(\mathcal A_{f}\mathcal N_q)$   for some prime $p$ and $q$ (see Subsection \ref{sect:pseudo} for  definitions). 

\subsection{The surface group conjecture.} Our last result represents some progress towards the solution of  Mel'nikov's conjecture. A finitely generated group $G$ is called {\bf Mel'nikov} if  all its subgroups of finite index are one-relator groups. The following problem appears in \cite[Problem 7.36]{Kou18} and also in \cite[Conjecture 2.17]{BauFine19}.
\begin{conj}[Mel'nikov] \label{Mel'conj} An infinite and  residually finite Mel'nikov group is either a free group, a surface group or a solvable Baumslag--Solitar group. 
\end{conj}
The two-generated case of the conjecture is proved by Gardam--Kielak--Logan \cite{GKL22}. In this paper we consider $n$-generator one-relator groups with $n\geq 3$. We say that a group is {\bf 2-free} if all its two-generated subgroups are free.

\begin{thmA} \label{Mel'Intro}
Let  $n\geq 3$ and let $F$ be the free group on $\{x_1,\ldots, x_n\}$. Let $  w\in F$ and suppose that $G=F/\langle\!\langle w\rangle\!\rangle$ is a residually finite   Mel'nikov group. Then  $G$ is 2-free. Moreover, if there exists a finite-index subgroup $H\leq G$ such that  $H^2(H; \Z)\neq 0$, then $G$ is a surface group.
\end{thmA}
This settles Conjecture \ref{Mel'conj} in the case of non-trivial second integral homology. Theorem \ref{Mel'Intro} also answers a question of Gardam--Kielak--Logan \cite[Question 3.3]{GKL22} on whether residually finite Mel'nikov groups have negative immersions (recall that $G$ has negative immersions if and only if $G$ is  2-free by \cite[Theorem 1.3 and Remark 1.7]{Lou22}).

\begin{rem*}
Marco Linton informed us that one can obtain the same conclusion in Theorem \ref{Mel'Intro} without assuming that $G$ is residually finite. This suggests that the Mel'nikov conjecture might also be true without assuming that $G$ is residually finite. This brings attention again to the two-generated case of the conjecture because in the proof of \cite{GKL22}, the condition that $G$ is residually finite is used in an essential way.
\end{rem*}

\subsection{More details on our methods} Before explaining the ideas that underpin the proofs, we make one observation about the assumptions of Theorem \ref{main}. 
Instead of   finite generability and vanishing of the second $L^2$-Betti number, the ingredients that we really  require in the proof  of Theorem  \ref{main}    are $G$ being $\FP_2(\F_p)$ and $p$-good at the variety of finite solvable groups (defined in Section \ref{goodsection}) to then conclude it from  Theorem \ref{main0}. However, it is in general very difficult to obtain prosolvable goodness, so we consider that the assumptions imposed in the statement of Theorem \ref{main} are not just sufficient, but also more natural and practical.  

{\bf From Theorem \ref{main0}  to Theorem \ref{main}:}. The   first author \cite[Theorem 1.1]{Ja22} showed that a group $G$ as in part (1) or part (2) of Theorem \ref{main} must be RFRS. So part (2)  implies part (1) and to deduce part (2) from Theorem \ref{main0} we only have to show that $G$ is prosolvable good and of type   $\FP_2(\F_p)$. Recall that if a  finitely presented  group $G$ is cohomologically good and has the same profinite completion as a surface group $S$ then, by L\"uck's approximation \cite{Luc94}, $b_2^{(2)}(G)=0$. Interestingly, in Claim \ref{criterion:cg}  we conversely prove that one can  recover goodness from $\cd(G)=2$ and $b_2^{(2)}(G)=0$.  This summarises how we prove  Theorem \ref{main} from Theorem \ref{main0}. 

{\bf The proof of Theorem \ref{main0}:} The ultimate goal is to prove that $G$ is Poincar\'e Duality of dimension two over the field $\F_p$ and then conclude from Bowditch \cite{Bow04}. During the proof we work with a smaller variety $\mathcal C$ than the variety of solvable groups to get a stronger formulation of Theorem \ref{main0}, but both equally  work and for simplicity we sketch the argument in terms of the solvable variety $\mathcal S.$

As we argue in Claim \ref{ends}, since $G_{\hat{\mathcal S}}$  is a freely indecomposable profinite group,   $G$ is also freely indecomposable and hence $H^1(G; \F_p[G])=0$. So it remains to show that $H^2(G; \F_p[G])\cong \F_p$ as trivial $G$-modules to prove that $G$ is Poincar\'e duality. There is a natural surjective map $H^2(G; \F_p[G])\lrar H^2(G; \F_p)\cong \F_p$ with kernel $K$. The remainder of the proof consists on understanding the cohomological, finiteness and profinite properties of $K$  to show that $K=0$. For example, using ideas from Sections \ref{cohomsection} and \ref{goodsection}, we can prove that the prosolvable completion $K_{\hat{\mathcal S}}$ of $K$ is zero. Later on, using some tools from Section \ref{L2section}, we deduce that $K$ is $L^2$-acyclic. Another fundamental ingredient is to ensure that $K$ has projective dimension at most two, i.e. that $\Ext_{\F_p[G]}^2(K, \F_p[G])=0$. For this, we compare $\Ext_{\F_p[G]}^2(K, \F_p[G])$ with $\Ext_{\Gsr}^2(K_{\hat{\mathcal S}}, \Gsr),$ which we already know to be zero because $K_{\hat{\mathcal S}}=0$.  Note that, in general, there is no way to make this comparison, since there is no canonical map between derived functors of $G$-modules and $G_{\hat{\mathcal S}}$-modules. Nevertheless, after deriving more information on the finiteness properties and on the cohomological goodness of the $G$-module $K$, we can prove that $\Ext_{\F_p[G]}^2(K, \F_p[G])\cong \Ext_{\Gsr}^2(K_{\hat{\mathcal S}}, \Gsr)$ using the tools from Section \ref{sect:pseudo}. Once we have proven the above, we know that there is a short exact sequence of right $\F_p[G]$-modules of the form \[0\lrar P_1\lrar P_0\lrar K\lrar 0,\]
where $P_0$ and $P_1$ are finitely generated projective $\F_p[G]$-modules with the same $\D_{\p[G]}$-dimension. Furthermore, since $K_{\hat{\mathcal S}}=0$, the map $P_1\lrar P_0$ is dense in the prosolvable topology. It follows from \cite[Theorem 4.3]{Ja22} that the map $P_1\lrar P_0$ is an isomorphism and so $K=0$. This completes the proof of Theorem \ref{main0}.

  {\bf The proof of Theorem \ref{surfaceword}:} This result is an  application of Theorem \ref{main} applied  to the double $G=F*_{u=u'}F'$  of $F$ along the corresponding word $u$.   

Lastly, we discuss the proof of Theorem \ref{Mel'Intro}.

{\bf Achieving 2-freeness:}  The first part of the proof of Theorem \ref{Mel'Intro} is ensuring that the one-relator group $G=F/\langle\!\langle w\rangle\!\rangle$   is 2-free. With this purpose, we  use a criterion of Louder--Wilton
  \cite[Theorem 1.5 and Definition 6.5]{Lou22}. This says that it suffices to rule out the existence of a subgroup $K$ of $F$ of rank two containing $w$ such that the canonical homomorphism $P=K/\langle\!\langle w\rangle\!\rangle \lrar G$ is an embedding of a non-free  2-generator one-relator group $P$ into $G$. We begin by  showing that $G$ is prosolvable good (which requires a different argument as the one in the proof of Theorem \ref{main}). This will lead to the closure $\overline{P}$ of $P$ in $\hat{G}$ being a projective profinite group. Then, by estimating the Betti numbers of the open subgroups of the closure $\overline{P}$ of $P$ in $\hat{G}$, we obtain that the $p$-Sylows of $\overline{P}$ are pro-$p$ cyclic. This implies that $\overline{P}$ is meta-abelian by a result of Zassenhaus \cite[Lemma 4.2.5]{Rib17} and hence that $P\cong B(1, m)$. A similar analysis on the virtual Betti numbers of $B(1, m)$ is needed to rule out this subgroup of $G$. 

{\bf Recognising the surface group:} The second part of Theorem \ref{Mel'Intro} consists on showing that $G$ is a surface group if there is a finite-index subgroup $H<G$ with $H^2(H; \Z)\neq 0$. For this we apply Theorem \ref{main0}, although it is not immediate to see that the required assumptions on $G$ are fulfilled. The group $G$ is a 2-free one-relator group and hence a virtually compact special group by Louder--Wilton \cite{Lou22} and Linton \cite{Lin22}. So $G$ is virtually $\mathcal{N}_p$-RFRS. It remains to prove that $G_{\hat{\mathcal S}}$ is Poincar\'e Duality of dimension two at some prime $p$. For this, we build on the fact that $G$ is prosolvable good. More precisely, we apply to the maximal $p$-quotients of $G_{\hat{\mathcal S}}$  the result of  Ando\v{z}ski\u{\i} \cite{An73} that a non-free one-relator pro-$p$ group that satisfies that all its open subgroups are one-relator groups must be a Demushkin group (recall that these are the Poincar\'e duality pro-$p$ groups of dimension two \cite[Section I.4.5, Example 2]{Ser97}). Interestingly, this was the original motivation of Mel'nikov to formulate Conjecture \ref{Mel'conj}.
 
\subsection*{Acknowledgments}

The work of the first author is partially supported by the grant     PID2020-114032GB-I00 of the Ministry of Science and Innovation of Spain  and by the ICMAT Severo
Ochoa project  CEX2019-000904-S4. The second author is supported by the Oxford--Cocker Graduate Scholarship. 

We wish to thank   Henry Wilton for  helpful conversations about Proposition \ref{double}. The second author would also like to thank Dawid Kielak and Richard Wade for their invaluable advice and their generous support.
\section{Preliminaries}
 
 \subsection{The Euler characteristic}
 All rings in this paper are associative and have the identity element. 
For a given ring  $R$, we will work with left  and right $R$-modules, but if we do not specify it, an $R$-module will mean a left $R$-module.

Let $G$ be a group and $k$ a commutative ring.
We denote by $k[G]$ the group ring of $G$ with coefficients in $k$ and by $I_{k[G]}$ its augmentation ideal.  
 
  Given  a ring $R$ and a $(\Z[G],R)$-bimodule $M$, we will denote the usual cohomology groups of $G$ with coefficients in $M$  by  $H^k(G; M)$. The group $H^k(G; M)$ is naturally a right $R$-module. 
  
  We say that an $R$-module $M$ is   $\FP$ if it has a finite resolution consisting of finitely generated   projective  $R$-modules. 
If $R\lrar \D$ is a division $R$-ring, the $\D$-Betti numbers of a left $R$-module $M$ are defined as 
\begin{equation} \label{Ddimension} b^{R,\D}_i(M)=\dim_{\D}  \Tor_i^{R}(\D,M).\end{equation}
 If $M$ is of type $\FP$, then  we put
\begin{equation}
\label{calculatingEuler}
\chi^{R, \D}(M)=\sum_{i} (-1)^i b^{R,\D}_i(M).
\end{equation}

  \begin{pro}\label{FP}
Let $0\lrar M_1\lrar M_2\lrar M_3\lrar 0$ be an exact sequence of $R$-modules. If any two of the modules $M_1$, $M_2$ or $M_3$ are $\FP$, then  the third one is also $\FP$. Moreover, if $\D$ is a division $R$-ring, then
$$\chi^{R, \D}(M_2)=\chi^{R, \D}(M_1)+\chi^{R, \D}(M_2).$$
\end{pro}
\begin{proof} 
 The first part of the proposition follows from
Schanuel's lemma \cite[Chapter VIII, Lemma 4.4]{Bro82}, and the second part is a consequence of the long exact sequence for $\Tor$ functors
 \cite[Chapter 2-3]{Weibel_HA}.  \end{proof}

Let $f\po R\lrar S$ be a ring homomorphism. We say that $f$ is  {\bf epic} if for every  ring $Q$ and homomorphisms $\alpha,\beta\po  S\lrar Q$, the equality $\alpha\circ f=\beta\circ f$ implies $\alpha=\beta$.

As in \cite[Chap. 7.2]{Co06}, we define an {\bf epic division $R$-ring} as a division ring $\mathcal{D}$ together with an epic homomorphism  $\varphi\colon R \lrar \mathcal{D}$. The condition on $\varphi$ to be epic is equivalent to the condition that the division closure of $\varphi(R)$ is equal to $\D$.

 Let $\varphi'\colon R \lrar \mathcal{D}'$ be  another epic division $R$-ring. A subhomomorphism of epic division $R$-rings $\mathcal D$ and $\mathcal D'$ is a homomorphism $\psi\colon K \lrar \mathcal{D}'$, where $K$ is a local subring of $\mathcal{D}$ containing $\varphi(R)$ with maximal ideal $\ker \psi$, such that $\psi\circ\varphi = \varphi'$. Two subhomomorphisms are equivalent if their restriction to the intersection of their domains coincide and are again subhomomorphisms. An {\bf specialization} $\mathcal{D} \lrar \mathcal{D}'$ of epic division $R$-rings is an equivalence class of subhomomorphisms in the previous sense. 
 \begin{pro}\label{equaleuler}
 Let $R$ be a ring,   $\varphi\colon R\lrar \D$ and $\varphi'\colon R\lrar\D^\prime$ two  epic division $R$-rings and $M$ a $\FP$-module. If there is a specialization $\mathcal{D} \lrar \mathcal{D}'$, then $\chi^{R, \D}(M)=\chi^{R, \D'}(M)$.
 \end{pro}
 \begin{proof}
 In view of Proposition \ref{FP}, it  is enough to prove the claim when $M=P$ is a finitely generated projective $R$-module, in which case $\chi^{R, \D_0}(M)=\dim_{\D_0} D_0\otimes_R P$ for any division $R$-ring $R\lrar \D_0$.

Let $\psi\colon K \lrar \mathcal{D}'$ be a homomorphism, where $K$ is a local subring of $\mathcal{D}$ containing $\varphi(R)$ with maximal ideal $\ker \psi$, such that $\psi\circ\varphi = \varphi'$.  Since projective modules over a local ring are free \cite{Kap58},  $K\otimes_R P\cong K^k$ as  $K$-modules  for some $k\in \N$.  Hence, $k=\dim_{\D} \D\otimes_R P$. Moreover, if we denote by $\mathcal Q\subset \D'$ the image of $\psi$ (which is the division quotient ring of $K$),  we have $\mathcal Q\otimes_R P\cong \mathcal Q\otimes_K K^k \cong \mathcal Q^k$ and hence $k=\dim_{\D'} \D'\otimes_R P$ as well.  
 \end{proof}

\subsection{Finiteness properties} \label{finiteSection}
  If a group has a finite classifying space, then it leads directly to a finite resolution of the trivial $\p[G]$-module $\p$ by finitely generated free $\p[G]$-modules. For our methods, it will be enough to work under the assumption of  being of type $\FP_2(\F_p)$, which we recall below.  
  
Let $R$ be a ring and let $n$ be either an integer or $n=\infty$. An $R$-module $M$ is {\bf  $\FP_n$}   if it admits a resolution by projective $R$-modules $P_*\lrar M\lrar 0$, where $P_i$ is finitely generated for all $0\leq i\leq n$. A group $G$ is  {\bf $\FP_n(\p)$} (for an integer $n\geq 0$ or for $n=\infty$)  if the trivial $\p[G]$-module $\p$ is   $\FP_n$. A group $G$ is said to be {\bf $\FP(\p)$} if the trivial $\p[G]$-module $\p$  admits a finite projective resolution.

We remind the reader that a group $G$ is $\FP_{\infty}(\p)$ if and only if it is $\FP_n(\p)$ for all $n$ by \cite[Chapter VIII, Proposition 4.5]{Bro82}; and also that a group $G$ of finite cohomological dimension is of type $\FP(\p)$ if and only if it is $\FP_{\infty}(\p)$ by \cite[Chapter VIII, Proposition 6.1]{Bro82}.
We also want to recall the following well-known result from \cite[Theorem 1]{Wat60} (see also \cite[Proposition 6.8]{Bro82}).

\begin{pro} \label{CohomologyA} Let $G$ be a $\FP(\p)$ group with $\cd(G)\leq n$; let  $R$ and $S$ be  two rings; and let $A$ be an  $(\p[G], R)$-bimodule and $B$ an $(R,S)$-bimodule. Then the canonical map 
$$ H^n(G; A)\otimes_R B\cong H^n(G; A\otimes_ R B).$$ is an isomorphism of right $S$-modules
\end{pro}
 
  We will often use the following consequence.
 
 \begin{cor}\label{corestriction}
 Let $G$ be a group of cohomological dimension two and let $U_1\leq U_2\leq G$ be subgroups of finite index in $G$. Then the corestriction map $H^2(U_1;\F_p)\lrar H^2(U_2;\F_p)$ is onto.
 \end{cor}

 \subsection{Pseudo-varieties of finite groups and profinite completions}\label{sect:pseudo}
We say that  a non-empty class of groups $\mathcal C$    is a  {\bf pseudovariety} if it is closed under subgroups, homomorphic images and finite direct products.  We denote by $\mathcal F$, $\mathcal P$, $\mathcal N$ and $\mathcal A$ 
the pseudovarieies  of finite groups,  polycyclic groups, finitely generated  nilpotent groups   and finitely generated abelian groups, respectively. If $\mathcal  C$ is a pseudovariety, we  put $\mathcal C_f=\mathcal C\cap \mathcal F$. The pseudovariety of all finite solvable groups $\mathcal P_f$ will be also denoted by  $\mathcal S$. If $p$ is a prime number,  then ${\mathcal C_p}$ will denote  the pseudovariety of finite $p$-groups lying in $\mathcal C$ and ${\mathcal C_{p^\prime}}$ will denote  the pseudovariety of finite $p^\prime$-groups lying in $\mathcal C$.  Given two pseudovarieties $\mathcal C $ and $\mathcal B$, we denote by $\mathcal C \mathcal B$ the pseudovariety consisting of isomorphic classes of group $G$ having a normal subgroup $N\in \mathcal C $ such that $G/N\in \mathcal B$.

Given a group $G$ and a pseudovariety $\mathcal C$ of finite groups we denote by $$G_{\widehat {\mathcal C}}=\varprojlim_{G/N\in \mathcal C}G/N$$ the pro-$\mathcal C$ completion of $G$.
 If $M$ is a  $\F_p[G]$-module, we write $L\le_{\mathcal C} M$ if $L$ is a submodule of $M$ of finite index and  $G/ St_G(M/L)\in \mathcal C$.
 The {\bf pro-$\mathcal C$ completion} $M_{\widehat{\mathcal C}}$ of $M$ is the inverse limit  
 $$M_{\widehat{\mathcal C}}=\varprojlim_{L\leq _{\mathcal C} M}M/L.$$ 
 If $M$ is a finitely generated $\F_p[G]$-module, then     we also have
 $$M_{\widehat{\mathcal C}}\cong \varprojlim_{G/N\in \mathcal C}M/\left (I_{\F_p[N]}M\right )\cong \varprojlim_{G/N\in \mathcal C} \left (\F_p[G/N]\otimes_{\F_p[G]} M \right ).$$ 
Observe that $M_{\widehat{\mathcal C}}$ becomes an $\F_p\llbracket G_{\widehat{\mathcal C}} \rrbracket $-module. We denote by $\delta_{M,\widehat{\mathcal C}}$ the canonical map 
\begin{equation} \label{notation:h2com}
     \begin{tikzcd}
          M\ar[r, "\delta_{M,\widehat{\mathcal C}}"] & M_{\widehat{\mathcal C}}. 
     \end{tikzcd}
\end{equation}   
We note that a right $\F_p[G]$-module has a pro-$\mathcal C$ completion in a similar way. 
If $M$ is a $\F_p[G]$-bimodule, then $M$ has  left and  right pro-$\mathcal C$ completions. In this case, to avoid ambiguity,  we  denote its right pro-$\mathcal C$-completion by $M_{\widehat{\mathcal C},r}$. 
The following result can be proved arguing as in \cite[Lemma 4.3]{GJ22}.

\begin{pro} \label{fpres}Let $M$ be a  $\F_p[G]$-module. Consider the natural  map $\alpha_{M}$ of $\F_p\llbracket G_{\widehat{\mathcal S}} \rrbracket $-modules 
\begin{equation*}
     \begin{tikzcd} 
           \F_p\llbracket G_{\widehat{\mathcal S}} \rrbracket  \otimes_{\F_p[G]} M\ar[r, "\alpha_{M}"] &  M_{\widehat{\mathcal C}}
     \end{tikzcd}
\end{equation*} defined by extending linearly the map $r\otimes m \mapsto r \cdot \delta_{M,\widehat{\mathcal C}}(m).$ Then the following holds. 
	\begin{enumerate}
		\item If $M$ is finitely generated, then $\alpha_{M }$ is onto.
		\item If $M$ is finitely presented, then $\alpha_{M }$ is an isomorphism.
		\item A  $\F_p[G]$-module morphism $\gamma\po M_1\lrar M_2$ induces the following commutative diagram of $ \Gcr$-modules.
\begin{equation*}
     \begin{tikzcd}
          \Gcr\otimes _{\F_p[G]} M_1 \ar[r, "\Id\otimes \gamma"] \ar[d, "\alpha_{M_1}"] &  \Gcr\otimes_{\F_p[G]} M_2 \ar[d, "\alpha_{M_2 }"]\\
         (M_1)_{\widehat{\mathcal C}}\ar[r, "\gamma_{\widehat{\mathcal C}}"] &  (M_2)_{\widehat{\mathcal C}}.
     \end{tikzcd}
\end{equation*}
  
	\end{enumerate}
	
\end{pro}

 Let $\Gamma$   be a finitely generated residually-$\mathcal C$ group.
 The {\bf $\mathcal C$-genus} of $\Gamma$, denoted by $\mathcal G_{\mathcal C}(\Gamma)$, is the set  of isomorphism classes of finitely generated residually-$\mathcal C$ groups $G$ having the same quotients in $\mathcal C$ as $\Gamma$. 
 It is  well-known (see, for example, \cite[Corollary 3.2.8]{RZ10}) that  $G\in \mathcal G_{\mathcal C}(\Gamma)$ if and only if  the pro-$\mathcal C$ completions of $G$ and $\Gamma$ are isomorphic: $G_{\widehat{\mathcal C}}\cong \Gamma_{\widehat{\mathcal C}}$.

In the next proposition, we should think of the inclusion $R\hookrightarrow S$ as the inclusion $\p[G]\hookrightarrow \Gcr$ for a group $G$ that is good in  a pseudovariety of finite groups $\mathcal C$.  
\begin{pro} \label{extmap} Let $R\hookrightarrow S$ be an extension of unital rings. Suppose that $A$ is an $R$-module such that $\Tor_i^R(S, A)=0$ for all $i\geq 1$ that admits a finite projective resolution of the form:
\begin{equation}\label{resol}
     \begin{tikzcd}[cramped, sep=small]  
    0\ar[r]&P_n\ar[r, "\partial_n"]& \cdots \ar[r, "\partial_1"]& P_0 \ar[r, "\partial_0"]& A\ar[r]& 0.
\end{tikzcd}
\end{equation}
Then there exists an isomorphism of right $S$-modules
\begin{equation*}
     \begin{tikzcd}[cramped, sep=small]  
     \Ext_R^n(A, R)\otimes_R S\ar[r]& \Ext_S^n(S\otimes_R A, S).
\end{tikzcd}
\end{equation*}
\end{pro}

To prove Proposition \ref{extmap}, we will need the following lemma, whose proof we leave as an exercise.
\begin{lem} \label{commutingfunctors} Let   $P$ be a projective $R$-module. Then there is a canonical isomorphism $F_P^S$ of $S$-modules
\begin{equation*}
    \begin{tikzcd}[cramped, sep=small]
        \Hom_R(P, R)\otimes_R S \ar[r, "F_P^S"]& \Hom_S(S\otimes_R P, S),
    \end{tikzcd}
\end{equation*}
that consists on extending linearly the map $\left(F_P^S(\psi\otimes s_1)\right)(s_2\otimes x)=s_2\cdot\psi(x)\cdot s_1$ for all $s_1, s_2\in S$, $x\in P$ and $\psi\in \Hom_R(P, R)$.  
\end{lem}

\begin{proof}[Proof of Proposition \ref{extmap}]

    After applying the functor $\Hom_R(-, R)$ to the  resolution (\ref{resol}), we get a chain complex of maps 
    $$\partial^i\po \Hom_R(P_{i-1}, R)\lrar \Hom_R(P_i, R)$$ and the abelian group $\Ext_R^n(A, R)$ is isomorphic to $\Hom_R(P_n, R)/\im \partial^n$. On the other hand, since $\Tor_i^R(S, A)=0$, we can apply the functor $(S\otimes_R -)$ to get a projective resolution of $S\otimes_R A$ of the form 
    \begin{equation*}
     \begin{tikzcd}[cramped, sep=small]  
    0\ar[r]&S\otimes_R P_n\ar[r, "\partial_n^S"]& \cdots  \ar[r, "\partial_1^S"]& S\otimes_R P_0 \ar[r, "\partial_0^S"]& S\otimes_R A\ar[r]& 0.
\end{tikzcd}
\end{equation*}
After applying the functor $\Hom_S(-, S)$ to this resolution, the induced chain complex has maps $\partial_S^i\po \Hom_S(S\otimes_R P_{i-1}, S)\lrar \Hom_S(S\otimes_R P_i, S)$. As before, the abelian group $\Ext_S^n(S\otimes_R A, S)$ is isomorphic to $\Hom_S(S\otimes_R P_n, S)/\im \partial_S^n$. The situation can be summarised with the following commutative diagram:
\begin{equation*}
     \begin{tikzcd}[cramped, sep=small]
          \Hom_R(P_{n-1}, R)\otimes_R S\ar[r] \ar[d, "F_{P_{n-1}}^S"] &  \Hom_R(P_n, R) \otimes_R S\ar[r] \ar[d, "F_{P_n}^S"] &  \Ext_R^n(A, R)\otimes_R S \ar[r] \ar[d] & 0\\
         \Hom_S(S\otimes_R P_{n-1}, S)\ar[r] &  \Hom_S(S\otimes_R P_n, S) \ar[r] & \Ext_S^n(S\otimes_R A, S) \ar[r] & 0.
     \end{tikzcd}
\end{equation*}
The maps $F_{P_i}^S$ are isomorphisms by Lemma \ref{commutingfunctors} and they naturally induce the required isomorphism  $\Ext_R^n(A, R)\otimes_R S\cong \Ext_S^n(S\otimes_R A, S).$
\end{proof}

\subsection{$L^2$-Betti numbers and their mod-$p$ variants for RFRS groups}  \label{L2section}
Let $G$ be a group. The usual Betti numbers of $G$ are defined by $$b_i(G)=b_i^{\Z[G],\Q}(\Z)= b_i^{\Q[G],\Q}(\Q)\textrm{\hspace{4pt} and \hspace{3pt} } b_{i,p}(G)=b_i^{\Z[G],\F_p}(\Z)=b_i^{\F_p[G],\F_p}(\F_p).$$

As part of our group cohomological tools, we require some estimations about first $L^2$-Betti numbers. They allow us to relate cohomology with rational and finite coefficients to the profinite completion of the group (for example, they are used in a criterion to recognise goodness  in Claim \ref{criterion:cg}). Since we are going  to consider only RFRS groups, we introduce the $L^2$-Betti numbers and the mod-$p$ $L^2$-Betti numbers of these groups in a convenient and uniform algebraic way. 

First recall the definition of RFRS groups. A group $G$ is called {\bf residually finite rationally solvable} or {\bf RFRS} if there exists a chain $G=H_0>H_1>\cdots$ of finite index normal subgroups of $G$ with trivial intersection such that $H_{i+1}$ contains a normal subgroup  $K_{i+1}$ of $H_i$ satisfying  that $H_i/K_{i+1}$ is    torsion-free abelian.  The chain $\{H_i\}$ satisfying this property is called   {\bf witnessing}. 

Let $\mathcal V$   an extension-closed pseudovariety of finite solvable groups. 
 We say that a group $G$ is  {\bf $\mathcal V$-RFRS} if 
 $G$ has a witnessing chain $\{H_i\}$ such that $H_i$ is normal in $G$ and $G/H_i\in \mathcal V$ for all $i$.   
  Observe that in this case $G$ is residually-$\mathcal V $.
 It is clear that  a RFRS group is also  $\mathcal S$-RFRS. We will pay special attention to the property of $\mathcal N_p$-RFRS, which is named $\mathrm{RFR}p$ in \cite{Kob20}.

Let $G$ be a RFRS group and $K$ a field, By \cite[Corollary 1.3]{And19}, we know that there exists a universal $K[G]$-ring of fractions $\D_{K[G]}$. The {\bf $L^2$-Betti} and {\bf mod-$p$ $L^2$-Betti numbers}
can be defined, in the sense of Equation (\ref{Ddimension}), as
$$b_i^{(2)}(G)=b_i^{\Q[G],\D_{\Q[G]}}(\Q) \textrm{\ \ and \ \ } b^{(2)}_{i,p}(G)=b_i^{\F_p[G],\D_{\F_p[G]}}(\F_p).$$
\begin{pro}\label{typeFP2}
Let $G$ be a RFRS group of cohomological dimension 2, let $K$ be a field and let $M$ be a $K[G]$-module of type $\FP$.

\begin{enumerate}
\item For any  two $K[G]$-division rings $\D_1$ and $\D_2$ we have  $\chi^{K[G], \D_1}(M)=\chi^{ K[G],\D_2}(M)$.
\item We have the inequality $$b_i^{K[G],\D_{K[G]}}(M)\leq  b_i^{K[G], K}(M).$$
\item If  $\beta_2^{(2)}(G)=0$, then  $G$ is of type $\FP_2(\Z)$.

\end{enumerate}
\end{pro}
 \begin{proof}
 (1) By \cite[Corollary 1.3]{And19}, $\D_{K[G]}$ is universal division $K[G]$-ring of fractions. Thus, there exists  specialization $\D_{K[G]}\lrar \D_i$ ($i=1,2$). Therefore, by Proposition \ref{equaleuler}, $\chi^{K[G], \D_1}(M)=\chi^{ K[G],\D_2}(M)$.
 
 (2) Again, this follows from the universality of $\D_{K[G]}$.
 
 (3) This follows from \cite[Corollary 3.4]{JL23}.
 \end{proof}
 In view of Part (1) of Proposition \ref{typeFP2}, if $G$ is a RFRS group, then we will simply  write $\chi^{K[G]}(M)$ instead of $\chi^{K[G], \D}(M)$.
 The {\bf Euler characteristic} $\chi(G)$ of a RFRS group $G$  of type $\FP$ is the Euler characteristic of the trivial $\F_p[G]$-module $\F_p$ (or, equivalently,  of the trivial $\Q[G]$-module $\Q$).

 The following theorem is the well-known L\"uck approximation. It shows, in particular, that $\b_1$ is a profinite invariant among finitely presented groups.
 \begin{thm}[\cite{Luc94}] \label{Luck} Let $G$ be residually finite $\FP_{n+1}(\Q)$ group and let $G=N_1>N_2>\dots >N_m>\dots $ be any sequence of finite-index normal subgroups  with $\bigcap_m N_m=1$. Then 
\[\lim_{m\rar \infty} \frac{b_n(N_m)}{|G:N_m|}= b_n^{(2)}(G).\]
\end{thm}

 The following result is a slight variation  of \cite[Theorem 4.3]{Ja22}.

\begin{thm} \label{Grothen} Let $\mathcal V $  be a   pseudovariety of finite solvable groups, $\mathcal E=\mathcal A\mathcal V$ and
  $G$ a $\mathcal V$-RFRS group.
Let $L$ and $M$ be  two  finitely generated right projective  $\p[G]$-modules   with $$ b_0^{\F_p[G],\D_{\F_p[G]}} (L)= b_0^{\F_p[G],\D_{\F_p[G]}}(M).$$
If  $f\po L\lrar M$ is a $\p[G]$-homomorphism such that the induced map 
\begin{equation*}
     \begin{tikzcd}[cramped, sep=small]
          L\otimes_{\p[G]} \Ger  \ar[r, "f_{\hat{\mathcal E}}"] & M\otimes_{\p[G]}\Ger
\end{tikzcd}
\end{equation*} 
is surjective, then $f$ is an isomorphism.
\end{thm}
\begin{proof}
There exists $L^\prime$ such that $L\oplus L^\prime$ is a free finitely generated $\F_p[G]$-module. Now apply    \cite[Theorem 4.3]{Ja22} to the map $L\oplus L^\prime\xrightarrow{f\oplus \Id} M\oplus L^\prime.$
\end{proof}
\subsection{Poincar\'e duality groups of dimension 2} 

By a {\bf surface group}   we will mean  the fundamental group of a compact closed surface of negative Euler characteristic. In the orientable case, surface groups admits the presentations of the form $$S_g=\langle x_1,\ldots,x_g, y_1,\ldots, y_g|\ [x_1y_1]\cdots [x_g,y_g]=1\rangle\  ( g\geq  2);$$ and in the  non-orientable closed case it is $$N_g=\langle x_1,\ldots, x_g|\ x_1^2\cdots x_g^2=1\rangle\ (g\geq 3).$$
Although free groups arise as fundamental groups of non-closed surfaces of negative Euler characteristic, we will not consider free groups as surface groups.

As we mentioned in the introduction, the main criterion we use to recognise surface groups is the fact that these are the Poincar\'e duality groups of dimension two as we recall below. To expand on this, we first recall this definition using Farrell's approach \cite{Far75}.
\begin{defi} \label{PD2def} Let $k$ be a commutative ring and let $G$ be a $\FP_2(k)$  group of cohomological dimension two. We say that $G$ is Poincar\'e Duality of dimension two over $k$ if $H^1(G; k[G])=0$ and $H^2(G; k[G])\cong k$ as  $k[G]$-modules.
\end{defi}

Hence the way we recognise if a group $G$ is a surface group from its prosolvable completion is by studying the coresponding cohomology groups $H^i(G, \p[G])$ to ensure that it is a Poincar\'e duality of dimension two (with $\p$ coefficients). In the case of integer coefficients, Eckmann--Linnell--M\"{u}ller \cite{Eck80, Eck83} showed that the group must be a surface group. Later on, Bowditch  \cite{Bow04} generalised this principle to arbitrary coefficients. 
\begin{thm}\label{PD2}
Suppose that $G$ be a Poincar\'e duality group of dimension two over the field $\F_p$. Then $G$ is a surface group.
\end{thm}
We have only stated the version over $\F_p$ because this is the one we use and it gives us stronger conclusions in Theorem \ref{surfaceword} than the analogous one on rational coefficients. However, our methods could be easily adapted to work over $\Q$ (by considering the adequate  skew-fields of Section \ref{L2section}). 

Let $F$ be a free group on $\{x_1,\ldots x_n\}$. We say that $w\in F$ is an {\bf orientable} (resp. {\bf non-orientable) surface word} if there exists an automorphism $\phi\in \Aut(F)$ such that $\phi(w)=[x_1,x_2][x_3,x_4]\cdots [x_{n-1}, x_n]$ (resp. $\phi(w)=x_1^2x_2^2\cdots x_n^2$). We shall also need the following standard result. For the convenience of the reader, we sketch its proof.

\begin{pro}\label{double}
Let $F$ be a finitely generated group of rank k and $u\in F$. Assume that $F*_{u=\overline{u}}\overline{F}$ is a surface group. Then $u$ is a surface word.
\end{pro}
\begin{proof} Suppose that there exists a closed surface $\Sg$ such that $F*_{u=\overline{u}}\overline{F}\cong \pi_1\Sg$. It is well-known that group-theoretic cyclic splittings  of $\pi_1\Sg$ must come from  a geometric splitting of the surface $\Sg$. This is a two-dimensional analogue of the so-called Stallings--Epstein--Waldhausen  construction for 3-manifolds  (see, for instance, the exposition of Shalen \cite[Section 2]{Sha01}). Hence   there exists an essential simple closed curve $\ga$ in $\Sg$ satisfying the following properties: 
\begin{itemize}
    \item The surface $\Sg_{\ga}$ that results from cutting $\Sg$ along $\ga$ has exactly two connected components $\Sg_1$ and $\Sg_2$ (and so the boundaries of both $\Sg_1$ and $\Sg_2$ have exactly one connected component). 
    \item There are isomorphisms $\psi_1\po \pi_1\Sg_1\lrar F$ and $\psi_2\po \pi_1\Sg_2\lrar \overline{F}$ such that $\psi_1(\partial\Sg_1)=u$ and $\psi_2(\partial\Sg_2)=\overline{u}$.
\end{itemize}
If $\Sg_1$ is a genus $g$ surface with one boundary component, then $\pi_1\Sg_1$ is  free in generators $a_1, \cdots, a_{k}$, where $k=2g$ (resp. $k=g+1$) if $\Sg$ is orientable (resp. non-orientable). Now we note that from  \cite[Theorem 1.5]{Cul81}  that the isomorphism of group pairs $(\pi_1\Sigma_1; \pi_1 \partial \Sigma_1)\cong (F; u)$ implies that the word $u\in F$ belongs to the   $\Aut(F)$-orbit of  the word $[a_1, a_2]\cdots [a_{k-1}, a_k]$ if $\Sg_1$ is orientable (resp. $a_1^2a_2^2\cdots a_k^2$ if $\Sg_1$ is non-orientable). Hence $u$ is a surface word.
\end{proof}

 \subsection{Cohomology of  profinite groups} \label{cohomsection}
Given a profinite group $\GG$ and a   discrete   $\p\llbracket \GG\rrbracket $-module $M$, we will denote the continuous cohomology groups (in the sense of Galois cohomology) of $\GG$ with coefficients in $M$ by $\Hc^k(\GG; M)$.


A {\bf topological $\GG$-module} $\MM$ is an abelian Hausdorff topological group $\MM$ which is endowed with the structure of an abstract $\GG$-module such that the action $\GG\times \MM\lrar \MM$ is continuous.  Given a profinite group $\GG$ and a topological  $\GG$-module $\MM$,  we denote by $\Hc^k(\GG; \MM)$ the $k$-th continuous cochain cohomology group of $\GG$ with coefficients in $\MM$ as in \cite[Definition 2.7.1]{Neu08}. For our computations, we are only interested in profinite modules. Recall that a {\bf profinite $\GG$-module} $\MM$ is a topological $\GG$-module which is additionally compact and totally disconnected. Equivalently, we can simply define a profinite $\GG$-module $\MM$ to be an inverse limit  of finite discrete $\GG$-modules. From this point of view, we can compute the cohomology of $\GG$ with coefficients in $\MM$ as follows. 

\begin{pro}\cite[Corollary 2.7.6]{Neu08} 
\label{inverselimit}
Let $\MM=\varprojlim_i M_i$ be an inverse limit of finite discrete $\GG$-modules $\{M_i\}_{i\in \N}$. Assume that $H^k(\GG;M_i)$ is finite for all
$i\in \N$, then the canonical map
\begin{equation*}
\begin{tikzcd}[cramped, sep=small]
    \Hc^{k+1}(\GG,\MM)\ar[r]&
     \varprojlim_i  \Hc^{k+1}(\GG, M_i)
\end{tikzcd}
\end{equation*} is an isomorphism.
\end{pro}

 \subsection{Cohomological goodness} \label{goodsection}
Let $G$ be an abstract group and let $\mathcal C$ be a pseudovariety of finite groups. We say that $G$ is  {\bf cohomologically $p$-good in $\mathcal C$} (or just $p$-good in $\mathcal C$) if for any finite discrete $ \p\llbracket G_{\widehat {\mathcal C}}\rrbracket$-module $M$ with $|M|=p^k$ and $G/St_G(M)\in \mathcal C$, the induced  map \begin{equation}\label{pgoodC}
\begin{tikzcd}[cramped, sep=small]
    \Hc^i(G_{\hat{\mathcal{C}}}; M)\ar[r]&
    H^i(G; M)
\end{tikzcd}
\end{equation} is an isomorphism. If, for all primes $p$, $G$ is $p$-good in  the variety of all finite groups, then we just say that $G$ is {\bf good} (recovering   Serre's original notion \cite{Ser97}).

\begin{pro} \label{tor} Let $\mathcal C$ be a pseudovariety of finite groups such that $\mathcal N_p\mathcal C=\mathcal C$. Let $G$ be a  $\FP(\F_p)$ group that is cohomologically good in $\mathcal C$, then $\Tor_{i}^{\p[G]}(\p, \Gcr)=0$ if $i\geq 1$. 
\end{pro}
\begin{proof}
The same proof as of \cite[Proposition 3.1]{Ja20} applies in this case.
\end{proof}
The previous proposition implies that, given a projective resolution of the trivial $\p[G]$-module $\p$ such as
\begin{equation*}
    \begin{tikzcd}[cramped, sep=small]
        0\ar[r] & P_n\ar[r]  &\cdots \ar[r]& P_0\ar[r] & \p\ar[r] & 0,
    \end{tikzcd}
\end{equation*}
we can apply the functor $(\Gcr\otimes_{\p[G]} -)$ to naturally obtain a projective resolution of the trivial $\Gsr$-module $\p$ of the form
\begin{equation*}
    \begin{tikzcd}[cramped, sep=small]
        0\ar[r] & \Gcr\otimes_{\p[G]} P_n\ar[r]  &\cdots \ar[r] & \Gcr\otimes_{\p[G]} P_0\ar[r] & \p\ar[r] & 0.
    \end{tikzcd}
\end{equation*}
Lastly, before moving on to the proofs of our results, we finish this section by recalling the definition of Poincar\'e Duality for profinite groups. Farrell's approach described in Definition \ref{PD2def}  naturally carries over to this setting.
\begin{defi} \label{PD2profinite}  Let $\GG$ be a profinite group of type $\FP_2(\F_p)$. We say that $\GG$ is a {\it Poincar\'e duality group of dimension $2$ over $\F_p$ at $p$} if $\cd_p(\GG)=2$, $H^1(\GG, \F_p \llbracket\GG\rrbracket)=0$ and $H^2(\GG; \F_p\llbracket \GG\rrbracket)\cong \F_p$  as abelian groups.    
\end{defi}

\section{The proof of Theorem \ref{main0}} \label{sect:rigidity} 
Theorem \ref{main0} can be strengthened to the following result. 
\begin{teo}\label{mainC0}
  Let $p$ be a prime, let $\mathcal V$ be an extension closed pseudovariety of finite solvable groups and let $\mathcal C=\mathcal N_p(\mathcal A_f \mathcal V)$. 
Let $G$  be a   $\mathcal V$-RFRS group     of type $\FP_2(\F_p)$  that  is $p$-good at $\mathcal C$ and has cohomological dimension 2. Suppose that $G_{\hat{\mathcal C}}$ is Poincar\'e Duality of dimension two at $p$. Then $G$ is a surface group.
 \end{teo}
 Theorem \ref{main0} corresponds to the case when  $\mathcal V$ (and hence $\mathcal C$) is the variety of finite solvable groups. Even if the condition of being $\mathcal V$-RFRS for a small $\mathcal V$ (such as $\mathcal N_p$) is stronger than that of being RFRS in the usual sense, the assumption on  $G_{\hat{\mathcal C}}$ being Poincar\'e Duality only involves quotients on $\mathcal C$. This  lead to stronger consequences not only on the profinite rigidity of surface groups (Theorem \ref{mainC}), but also on the topic of measure equivalent words  (Theorem \ref{surfacewordstrong}).
 \begin{proof}[Proof of Theorem \ref{mainC0}]
  By Kerckhoff's realisation theorem \cite{Ker83}, a  torsion-free virtually surface group is a surface group itself. So we can change $G$ by a subgroup of finite index and assume that $H^2(G_{\widehat{\mathcal C}}; \F_p\llbracket G_{\widehat{\mathcal C}}\rrbracket)\cong \F_p$ is a trivial $G_{\widehat{\mathcal C}}$-module. 
 
Our strategy is to show that $G$ is $\PD_2$ over $\F_p$ and conclude from Theorem \ref{PD2} that $G$ is a surface group. According to Definition \ref{PD2def}, this amounts to check that $H^1(G; \p[G])=0$ and that $H^2(G; \p[G])\cong \p$.
  \begin{claim}\label{ends} The group $G$ satisfies $H^1(G; \p[G])=0$.
\end{claim}
\begin{proof}
 Notice that $G_{\widehat {\mathcal C}}$ does not split as a free pro-$\mathcal C$ product because, otherwise, if we had $G_{\widehat {\mathcal C}}=\UU\coprod \VV$ for non-trivial $\UU$ and $\VV$, it would follow that $\cd_p(\UU)\leq 1$ and $\cd_p(\VV)\leq 1$  by a result of Serre \cite[Section 1.4]{Ser97}. This would contradict the fact that $\cd_p(G_{\widehat {\mathcal C}})=2$. So both $G_{\widehat {\mathcal C}}$ and $G$ are freely indecomposable. In particular, since $G$ is not virtually cyclic, we deduce from Stallings theorem \cite{St68} that $G$ is one-ended. Therefore,  $H^1(G; \p[G])$ vanishes  (see, for example, the book of Dicks--Dunwoody \cite[Theorem 6.10]{Dic89}).  \renewcommand\qedsymbol{$\diamond$}
\end{proof}
 Since $G$ if of type $\FP_2(\F_p)$,  we have an exact sequence.
 
 \begin{equation}\label{eq1}
 \begin{tikzcd}[cramped, sep=small]
    0\ar[r]&  P_2\ar[r, "\alpha"] & P_1\ar[r, "\beta"] &  P_0 \ar[r]& \p\ar[r]& 0,
\end{tikzcd}
\end{equation}
where the  $P_i$ are finitely generated projective $\p[G]$-modules.
We name $$M:=H^2(G;\p[G]).$$ After applying the right-exact contravariant functor $\Hom(-, \p[G])$, we obtain   the following exact sequence of  right $\p[G]$-modules:
\begin{equation}\label{e3}
\begin{tikzcd}[cramped, sep=small]
    0\ar[r]&  Q_0\ar[r, "\beta^*"] &Q_1 \ar[r, "\alpha^*"] & Q_2\ar[r]& M\ar[r]& 0,
\end{tikzcd}
\end{equation}
where $Q_i=\Hom(P_i, \p[G])$.
It is exact by  Claim \ref{ends}.  In particular, $M$ is $\FP$ (as a right $\p[G]$-module). 

\begin{claim} \label{h2com} Let $A$ be a  $\p[G]$-bimodule that is finitely generated  as  a right  $\p[G]$-module.  Then  $H^2(G; A)_{\widehat {\mathcal C}}$ and $ \Hcp^2(\Gc;  A_{\widehat{\mathcal C},r})$ are isomorphic as   right  $\p[G]$-modules (the notation for $A_{\widehat{\mathcal C},r}$ was introduced after (\ref{notation:h2com})).
\end{claim}
\begin{proof} First observe that since $A$ is finitely generated as right $\F_p[G]$-module, the right   $\F_p[G]$-module $H^2(G;A)$,  which is isomorphic by Proposition  \ref{CohomologyA} to $M\otimes_{\F_p[G]} A$ and so is also finitely generated. Hence 
$$A_{\widehat{\mathcal C},r}\cong \varprojlim_{G/N\in \hspace{1pt} \mathcal C} A \otimes_{\F_p[G]}\F_p[G/N]$$ and $$
H^2(G; A)_{\widehat {\mathcal C}}\cong \varprojlim_{G/N\in  \hspace{1pt} \mathcal C} H ^2(G; A) \otimes_{\F_p[G]}\F_p[G/N].  $$

 Again using that $A$ is finitely generated as a right $\F_p[G]$-module, we obtain that $A \otimes_{\F_p[G]}\F_p[G/N]$ is finite if $N$ is a normal subgroup of $G$ of finite index. 
 Thus, $\Hc^1 (G_{\widehat{\mathcal C}},A \otimes_{\F_p[G]}\F_p[G/N])$ is finite. 
 Hence, 
 by Proposition \ref{inverselimit}, the canonical map 

 \begin{equation} \label{4.3eq1}
\begin{tikzcd}[cramped, sep=small]
     \Hc^2 (G_{\widehat{\mathcal C}};A _{\mathcal C,r})\ar[r]&
     \varprojlim_{G/N\in \mathcal C} \Hc^2 (G_{\widehat{\mathcal C}};A \otimes_{\F_p[G]}\F_p[G/N])
\end{tikzcd}
\end{equation} 
 is an isomorphism.

From the $p$-goodness of $G$ in $\mathcal C$ we get that  the restriction map 
\begin{equation} \label{4.3eq2}
\begin{tikzcd}[cramped, sep=small]
      \Hc^2 (G_{\widehat{\mathcal C}};A \otimes_{\F_p[G]}\F_p[G/N])\ar[r]&
    H^2 (G ;A \otimes_{\F_p[G]}\F_p[G/N])
\end{tikzcd}
\end{equation} 
  is an isomorphism.

Finally, by Proposition \ref{CohomologyA}, the canonical map \begin{equation} \label{4.3eq3}
\begin{tikzcd}[cramped, sep=small]
      H^2(G;A) \otimes_{\F_p[G]}\F_p[G/N]\ar[r]&
     H^2(G;A\otimes_{\F_p[G]}\F_p[G/N])
\end{tikzcd}
\end{equation}  is also an isomorphism. The claim follows from (\ref{4.3eq1}), (\ref{4.3eq2}) and (\ref{4.3eq3}). \renewcommand\qedsymbol{$\diamond$}
 \end{proof}

We can compute  the pro-$\mathcal C$ completion of $M$ directly from Claim \ref{h2com}.

\begin{claim}\label{Mcom}
We have that $  M_{\widehat{\mathcal C}} \cong \p$, as trivial right $G$-modules. 
\end{claim}
\begin{proof}
This follows from  Claim \ref{h2com}, because $\Hc^2(G_{\hat{\mathcal{C}}};\p \llbracket G_{\widehat {\mathcal C}}\rrbracket )\cong \p$. \renewcommand\qedsymbol{$\diamond$}\end{proof}

By the $p$-goodness of $G$ at $\mathcal C$, we know that $H^2(G,\F_p)\cong \Hc^2(\Gc,\F_p)\cong \F_p$.
Consider the kernel $K$ of the natural surjective map  $H^2(G; \p[G])\lrar H^2(G; \p)$.  Observe that, since $M$ and $\F_p$ are $\FP$,  $K$ is $\FP$ by Proposition \ref{FP}.

\begin{claim} \label{Kcomp} The pro-$\mathcal C$ completion $K_{\widehat{\mathcal C}} $ is zero.
\end{claim}
\begin{proof} 

 Let $L$ be a submodule of $K$ of finite index such that $G/St_G(K/L)\in \mathcal C$. Observe that
 $M/K\cong H^2(G; \p)\cong  \F_p.$ Hence  $G/St_G(M/L)\in \mathcal C$ because $\mathcal {N}_p \mathcal C=\mathcal C$. 
Therefore, the canonical map $K_{\widehat{\mathcal C}} \lrar M_{\widehat{\mathcal C}} $ is injective. By Claim \ref{Mcom}, the canonical map $M_{\widehat{\mathcal C}}\lrar\F_p$ is an isomorphism. Hence   $K_{\widehat{\mathcal C}} $ should be zero. \renewcommand\qedsymbol{$\diamond$}
 \end{proof}

This is the first of a sequence of claims where we aim to show that the module $K$ is zero. Since $G$ is cohomologically $p$-good in $\mathcal{C}$, by Proposition \ref{tor}, $$\Tor_i^{\p[G]}(\Gcr, \p)=0 \textrm{ for all $i\geq 1$}.$$ Hence we can apply the functor $(\Gcr\otimes -)$ to the $\p[G]$-resolution of the trivial module $\p$ of equation (\ref{eq1}) to get the coresponding projective resolution of the trivial $\Gcr$-module $\p$:
 
\begin{equation} \label{e4}
     \begin{tikzcd}[cramped, sep=small]
        0\ar[r] &  (P_2)_{\widehat{\mathcal C}} \ar[r, "\hat \alpha"]&(P_1)_{\widehat{\mathcal C}} \ar[r, "\hat \beta"] &  (P_0)_{\widehat{\mathcal C}}\ar[r]& \p\ar[r]& 0.
     \end{tikzcd}
\end{equation}

\begin{claim}\label{TorM} We have $\Tor_i^{\p[G]}(M, \Gcr)=0$ for all $i\geq 1$.
\end{claim}

\begin{proof}
Since $\Gc\cong S_{\widehat {\mathcal C}}$,  \[\Ext_{\Gcr}^2 (\p, \Gcr)\cong \Hc^2(\Gc; \p)\cong {\p}.\] 
We can  use the projective resolution $C_*$ of $M$ given by (\ref{e3}) to compute the group $\Ext_{\p[G]}^2 (M, \p[G])$. 
It  will be isomorphic to
 $$H^2\left(\Hom_{\p[G]}(C_*, \p[G])\right) \cong P_0/\im \beta\cong \p.$$ 
Similarly, we will compute $\Ext_{\Gcr}^2 (M_{\hat{\mathcal{C}}}, \Gcr).$ 
For this purpose, recall the isomorphisms  $\Hc^0(\Gc, \Gcr)=0$, $\Hc^1(\Gc, \Gcr)=0$ and $\Hc^2(\Gc, \Gcr)\cong M_{\hat{\mathcal{C}}}\cong \F_p$.

The equation (\ref{e4}) is a projective resolution of $\p$, so after applying the functor $\Hom_{\Gcr}(-, \Gcr)$ to (\ref{e4}), 
we get the following exact sequence of right $ \Gcr$-modules. 
 
\begin{equation}\label{e6}
\begin{tikzcd}[cramped, sep=small]
    0\ar[r]& (Q_0)_{\widehat{\mathcal C}}\ar[r, "\beta_{\hat{\mathcal C}}^*"] &
    (Q_1)_{\widehat{\mathcal C}}
    \ar[r, "\alpha_{\hat{\mathcal C}}^*"] & (Q_2)_{\widehat{\mathcal C}}\ar[r]& M_{\hat{\mathcal C}}\ar[r]& 0.
\end{tikzcd}
\end{equation}

By Lemma \ref{commutingfunctors}, the previous resolution is what we would get if we apply first the functor $\Hom(-, \p[G])$, obtaining the projective resolution of $M$ described in (\ref{e3}), and then we apply $-\otimes_{\p[G]} \Gcr $. This implies the claim.  \renewcommand\qedsymbol{$\diamond$}
\end{proof}

\begin{claim} \label{extK}
We have that $\Ext_{\p[G]}^2 (K, \p[G])=0$.
\end{claim}
\begin{proof} 
Since $\cd(G)=2$,  when applying the derived functor $\Ext(-, \p[G])$ to the short exact sequence of $\p[G]$-modules $0\lrar K\lrar  M\lrar \p\lrar 0$ we get an exact sequence of right $\p[G]$-modules as follows.
\begin{equation*}
     \begin{tikzcd}[cramped, sep=small]
    \Ext_{\p[G]}^2 (\p, \p[G])\ar[r]&\Ext_{\p[G]}^2 (M, \p[G])\ar[r]  & \Ext_{\p[G]}^2 (K, \p[G])\ar[r] & 0.
\end{tikzcd}
\end{equation*}
We showed that $\Ext_{\p[G]}^2 (M, \p[G])$ is isomorphic to the trivial $\p[G]$-module $\p$. 
In particular, its quotient $\Ext_{\p[G]}^2 (K, \p[G])$ is a finite $\p[G]$-module with a trivial action. This immediately gives the isomorphism
\begin{equation} \label{extK1}
 \Ext_{\p[G]}^2(K, \p[G])\cong \Ext_{\p[G]}^2(K, \p[G])\otimes_{\p[G]} \Gcr.\end{equation} 
By the long exact sequence of $\p[G]$-modules that we obtain after applying the functor  $-\otimes_{\F_p[G]}\Gcr$ 
to the short exact sequence $0\lrar K\lrar M\lrar \p\lrar 0$, we get that $\Tor_i^{\p[G]}(\Gcr, K)=0$
 (since this holds for the modules $\p$ and $M$ by the pro-$\mathcal C$ $p$-goodness of $G$ and Claim \ref{TorM}). 
 By Proposition \ref{extmap},
\begin{equation} \label{extK2}  \Ext_{\p[G]}^2(K, \p[G])\otimes_{\p[G]} \Gcr\cong\Ext_{\Gcr}^2(K_{\hat{\mathcal{ C}}}, \Gcr).\end{equation}
The previous is zero because  $K_{\hat{\mathcal{ C}}}=0$ by  Claim \ref{Kcomp}. Thus,  $\Ext_{\p[G]}^2 (K, \p[G])=0$ by the isomorphisms (\ref{extK1}) and (\ref{extK2}).  \renewcommand\qedsymbol{$\diamond$}
\end{proof}

\begin{claim} \label{chiM} The $\p[G]$-modules  $\p$ and $M$ satisfy that $$\chi^{\F_p[G]}(\p)=\chi^{\F_p[G]}(M).$$
\end{claim}
\begin{proof} This follows from  the equations (\ref{eq1}) and (\ref{e3}).  \renewcommand\qedsymbol{$\diamond$}
\end{proof}

\begin{claim} \label{L2K} We have that $b_i^{\F_p[G],\D_{\F_p[G]}}(K)=0$  for all $i\geq 0$.
\end{claim}
\begin{proof} On the one hand, Claim \ref{Kcomp}  implies that $K\otimes_{\p[G]}\F_p=0$. By Proposition \ref{typeFP2}(2),  
 $$b_i^{\F_p[G],\D_{\F_p[G]}}(K)\leq  b_i^{\F_p[G], \F_p}(K)$$ and hence $b_0^{\F_p[G],\Dp}(K) =0$. 

On the other hand, from the additivity of $\chi^{\F_p[G]}$ (Proposition \ref{FP}) and the short exact sequence 
\begin{equation*}
     \begin{tikzcd}[cramped, sep=small]
    0\ar[r]& K\ar[r] & M\ar[r] & \p\ar[r] & 0,
\end{tikzcd}
\end{equation*}
  we get that $ \chi^{\F_p[G]}(K)-\chi^{\F_p[G]}(M)+\chi^{\F_p[G]}(\p)=0$. Moreover,  $\chi^{\F_p[G]}(\p)=\chi^{\F_p[G]}(M)$ by Claim \ref{L2K}, and so $ \chi^{\F_p[G]}(K)=0$. 
  
  Since $K$ has projective dimension at most two by Claim \ref{extK},  we have $ b_i^{\F_p[G],\D_{\F_p[G]}}(K)=0$ for all $i\geq 2$. In addition, by formula (\ref{calculatingEuler}), $$0=\chi^{\p[G]}(K)=-b_1^{\F_p[G],\D_{\F_p[G]}}(K) + b_0^{\F_p[G],\D_{\F_p[G]}}(K). $$ So it also follows that $b_1^{\F_p[G],\D_{\F_p[G]}}(K) =0$. \renewcommand\qedsymbol{$\diamond$}
\end{proof}

\begin{claim} \label{Ktrivial} The $\p[G]$-module $K$ is zero.
\end{claim}
\begin{proof}
Since $ \Ext_{\p[G]}^2 (K, \p[G])= 0$, the finitely generated right $\p[G]$-module $K$ has projective dimension at most two. Hence we have a short exact sequence
\begin{equation*}
     \begin{tikzcd}[cramped, sep=small]
    0\ar[r]& T_1\ar[r, "\gamma"] & T_0\ar[r] & K\ar[r] & 0,
\end{tikzcd}
\end{equation*} where $T_0$ and $T_1$ are finitely generated  projective $\p[G]$-modules. We know that $ b_0^{\F_p[G],\D_{\F_p[G]}}(K) =b_1^{\F_p[G],\D_{\F_p[G]}}(K) =0$ by Claim \ref{L2K}, so it follows that the induced map 
\begin{equation*}
     \begin{tikzcd}[cramped, sep=small]
      T_1  \otimes_{\p[G]}\mathcal D_{\p[G]}\ar[r] &  T_0  \otimes_{\p[G]}\mathcal D_{\p[G]}
\end{tikzcd}
\end{equation*} 
is an isomorphism. Hence $ b_0^{\F_p[G],\D_{\F_p[G]}} (T_1)= b_0^{\F_p[G],\D_{\F_p[G]}}(T_0)$. Furthermore, $K_{\widehat{\mathcal C}}=0$ by Claim \ref{Kcomp}, so the induced map \begin{equation*}
     \begin{tikzcd}[cramped, sep=small]
     T_1 \otimes_{\p[G]}\Gcr\ar[r] &  T_0 \otimes_{\p[G]}\Gcr
\end{tikzcd}
\end{equation*} 
 is surjective. Applying Theorem \ref{Grothen},  we obtain that $\gamma$ is an isomorphism. Thus $K$ is trivial.  \renewcommand\qedsymbol{$\diamond$}\end{proof}

We have shown that the kernel $K$ of the natural map $H^2(G; \p[G])\lrar H^2(G; \p)$ is zero and hence $H^2(G; \p[G])\cong \p$. 
By Theorem \ref{PD2}, this implies that $G$ is a surface group.  \end{proof}

\section{The proof of Theorem \ref{main}}
   It was proven in \cite[Theorem 1.1 and Proposition 4.1]{Ja22} that a group $G$ satisfying the assumptions of Theorem \ref{main} is   $\mathcal N_q$-RFRS and residually-$\mathcal N_p$ for all primes $p$ and $q$   if $S$ is an orientable surface group and $q>2$ if  $S$ is a  non-orientable surface group. So Theorem \ref{main} is implied by the following stronger result. 

 \begin{teo}\label{mainC}
 Assume that either
\begin{enumerate}
\item $S$ is an orientable surface group  and $q$ and $p$ are  arbitrary primes or
\item $S$ is a non-orientable surface group,  $p=2$ and $q\neq 2$ is prime.
\end{enumerate}
  Let   $G$   be a $\mathcal N_q$-RFRS group and assume  also that $G$ is  residually-$\mathcal N_p$,  finitely generated, of cohomological dimension 2 and has trivial second $L^2$-Betti number. Let $\mathcal C=\mathcal N_p(\mathcal A_f \mathcal N_q)$. If $G_{\widehat{\mathcal C}}\cong S_{\widehat{\mathcal C}}$, then $G\cong S$.
 \end{teo}

  We want to use  Theorem \ref{mainC0}.    First observe  that $G$ is of type $\FP_2(\Z)$ by Proposition \ref{typeFP2}(3). The next main ingredient is the following.

\begin{claim}\label{criterion:cg}  The group $G$ is cohomologically $p$-good in $\mathcal C$.
\end{claim}
 
\begin{proof} 

Weigel and Zalesskii \cite[Proposition 3.1]{Wei04} showed that an abstract group $G$ all of whose finite-index subgroups $H<G$ are {\bf cohomologically pro-$p$ good} (i.e. the canonical map $\Hc^i(H_{\hp}; \p)\lrar H^i(H; \p)$ is an isomorphism for all integers $i\geq 0$) must be cohomologically $p$-good. In our setting,  the same argument works to show that $G$ is cohomologically $p$-good in $\mathcal C$ if, for all $i$ and for all finite-index subgroups $H<G$ that are open in the pro-$\mathcal C$ topology, the natural map $\Hc^i(H_{\hp}; \p)\lrar H^i(H; \p)$ is an isomorphism.  

The previous map is always an isomorphism for $i\leq 1$ and, under our assumptions on  $G$ and $\Gc$, both $\Hc^i(H_{\hp}; \p)$ and $H^i(H; \p)$ vanish if $i\geq 3$.  Furthermore, the natural map $\Hc^2(H_{\hp}; \p)\lrar H^2(H; \p)$ is always injective by \cite[Section 2.6]{Ser97}. So it suffices to show that $\Hc^2(H_{\hp}; \p)\cong H^2(H; \p)$. 


There exists a surface group $S^\prime$ such that $\Hcc\cong S^\prime_{\hat{\mathcal N_p}}$. In particular, $b_{1, p}(H)=b_{1,p}(S^\prime)$. 
By Theorem \ref{Luck}, $\b_1(H)=\b_1(S^\prime),$ and so (\ref{calculatingEuler}) implies that $\chi (H)=\chi(S^\prime)$.

Thus, taking again into  account the formula (\ref{calculatingEuler}), we obtain 
\begin{multline*}
b_{2, p}(H)=\chi(H)+b_{1,p}(H)-b_{0,p}(H)= \chi(S^\prime)+b_{1,p}(S^\prime)-b_{0,p}(S^\prime)=\\
b_{2,p}(S^\prime)=\dim_{\p}\Hc^2(S^\prime_{\hp}; \p)=\dim_{\p}\Hc^2(H_{\hp}; \p).\end{multline*}
This  shows that the natural injective map $\Hc^2(H_{\hp}; \p)\lrar H^2(H; \p)$ is an isomorphism. \qedhere\end{proof}
Since $S$ is pro-$\mathcal C$ good at $p$,  its pro-$\mathcal C$ completion $S_{\hat{\mathcal C}}$ is Poincar\'e duality of dimension two at $p$. So $G$ lies under the assumptions of Theorem \ref{mainC0} and hence $G$ is a surface group. Surface groups are distinguished from each other by their abelianisation. Thus $G\cong S$ and the proof of Theorem \ref{mainC} is complete.

\section{The proof of Theorem \ref{surfaceword}}
Theorem \ref{surfaceword} is a consequence of the following stronger result.
 \begin{teo} \label{surfacewordstrong}
Assume that either
\begin{enumerate}
\item $F$ is a free group freely generated by generators $x_1,y_1,\ldots, x_d,y_d$,   $w=[x_1,y_1]\cdots [x_d,y_d]$ and $q$ and $p$ are  arbitrary primes or
\item $F$ is a free group freely generated by generators $x_1,\ldots, x_k$,  $w=x_1^2\cdots x_k^2$,  $p=2$ and $q$ is a prime different from 2.
\end{enumerate}
 Let  $u\in F$. Assume that $w$ and $u$ are measure equivalent in the pseudovariety ${\mathcal N_{p}(\mathcal A_{f}}\mathcal N_q)$.  Then there exists $\phi\in \Aut(F)$ such that  $u=\phi(w)$.

\end{teo}

    \begin{proof}[Proof of Theorem \ref{surfacewordstrong}] We want to apply Theorem \ref{mainC} and the following two claims allow us to do so.
 \begin{claim}
The element $u$ is not a proper power in $F$.
\end{claim}

\begin{proof}
Assume that $u$ is an $r$-power for some prime $r$.
Define $\mathcal U=\mathcal A_r\mathcal N_q$. Then by \cite[Proposition 3.2]{GJ22}, $w$ is   a  $r$-power  in $\mathbf G= F_{\widehat {\mathcal U}}$. It is clear then $r\ne q$. Thus, we assume that $r\ne q$.

Observe that since $w$ is an $r$-power,  it is a pro-$q$ element  of $\mathbf G$.  But this  is not the case. Let us prove it, for example, if $w=[x_1,y_1]\cdots [x_d,y_d]$. In this case consider the normal subgroup $\mathbf N$ of $\mathbf G$ generated by $x_1^q, y_1,x_2,\ldots, y_d$. Then $w[\mathbf N,\mathbf N]=[x_1,y_1][\mathbf N,\mathbf N]$ is not a pro-$q$ element in $\mathbf G/[\mathbf N,\mathbf N]$. \renewcommand\qedsymbol{$\diamond$}
\end{proof}
Consider the doubles $G=F*_uF$ and $S=F*_wF$. The group $G$ is an one-relator group without torsion. So $G$ has cohomological dimension 2 by \cite{Lyndon50} and trivial second $L^2$-Betti number by \cite{Dic07}.
In addition, by \cite{Ba62}, $G$ is residually-$\mathcal N_r$  for any prime $r$. Since $$G_{\widehat{\mathcal N_q}}\cong (F*_wF)_{\widehat{\mathcal N_q}},$$
arguing as in the proof of  \cite[Proposition 4.1]{Ja22} we obtain that $G$ is $\mathcal N_q$-RFRS. 

Let $\mathcal C= {\mathcal N_{p}(\mathcal A_{f}}\mathcal N_q)$. 

\begin{claim} There is an isomorphism $G_ {\widehat {\mathcal C}}\cong S_{\widehat {\mathcal C}}$.
\end{claim}

\begin{proof}
In view of \cite[Corollary 3.2.8]{RZ10}, we have to show that $G$ and $S$ have  the same finite quotients in $\mathcal C$. Let $P\in \mathcal C$. Denote  by $h(G,P)$ (resp. $e(G,P)$) the number of  homomorphims  (resp. epimorphisms) $H\lrar P$. Then we have 
$$h(G,P)= \sum_{a\in P} |u_P^{-1}(a)|^2=\sum_{a\in P} |w_P^{-1}(a)|^2=h(S,P).$$ Taking into account that
$$h(G,P)=\sum_{T\leq P} e(G,T) \textrm{\ and \ } h(S,P)=\sum_{T\leq P} e(S,T),$$ and arguing by induction on $|P|$, we obtain that $e(G,P)=e(S;P)$ for every $P\in \mathcal C$. Thus, $G$ and $S$ have  the same finite quotients in $\mathcal C$.\renewcommand\qedsymbol{$\diamond$}
 \end{proof}
    By Theorem \ref{mainC}, we obtain that $G\cong S$.
 Thus, by Proposition \ref{double}, $u$ is a surface word.
\end{proof}

 \section{Mel'nikov's groups}
 In this section we prove Theorem \ref{Mel'Intro} from the introduction, but we first restate a more complete formulation.
 \begin{thm} \label{Mel'2}
     Let  $n\geq 3$ and $F$ the free group on $\{x_1,\ldots, x_n\}$. Let $  w\in F$ and let $G=F/\langle\!\langle w\rangle\!\rangle$ be a residually finite   Mel'nikov group. Then  $G$ is 2-free. Moreover, the following statements hold.
     \begin{itemize}
         \item[(I)] If   $H^2(H; \F_p)=0$ for all finite-index subgroups $H\leq G$ and all primes $p$, then $\widehat G$ is a projective profinite group.  
         \item[(II)] Otherwise, if $H^2(H; \F_p)\neq 0$ for some finite-index subgroup $H\leq G$ and some prime $p$, then $G$ is a surface group.
     \end{itemize}
 \end{thm}
 \begin{proof}[Proof of Theorem \ref{Mel'2}]
 We divide the proof   in various claims. We fix a prime $p$ for the rest of the discussion. Without loss of generality we will assume that $w\ne 1$ and $w$ is not primitive. It is clear also that $w$ is not a proper power. Lyndon's asphericity theorem \cite{Lyndon50} proves that the presentation 2-complex of $G=F/\langle\!\langle w\rangle\!\rangle$ is a $K(G, 1)$ and that $\chi(G)=2-n$.

 \begin{claim}\label{abfi}
 For any subgroup $H\leq G$ of finite index we have  $$\dim_{\F_p} H^1(H;\F_p)\leq (n-2)|G:H|+2.$$
 \end{claim}

 \begin{proof} Recall that $\chi(H)=\dim_{\F_p} H^2(H;\F_p)-\dim_{\F_p} H^1(H;\F_p)+1$ for any finite-index subgroup $H\leq G$. Moreover, it is clear that   $\dim_{\F_p} H^2(H;\F_p)\leq 1$ (see \cite{Lyndon50}). The multiplicativity of the Euler characteristic  $\chi(H)=(2-n)|G:H|$ gives the desired conclusion. \renewcommand\qedsymbol{$\diamond$}
 \end{proof}

  Our next aim consists on showing that $G$ is $p$-good in $\mathcal C$. Contrary to the assumptions of Claim \ref{criterion:cg}, we do not have any strong input on the pro-$\mathcal C$ completion of $G$, so we argue slightly differently.  
 
 \begin{claim} \label{restriction}
 Let $U_1\leq U_2\leq G$ be two subgroups of finite index in $G$, such that $p$ divides $|U_2:U_1|$. Then the restriction map $H^2(U_2, \F_p)\lrar H^2(U_1,\F_p)$ is trivial.
 \end{claim}
 \begin{proof}  
 If one of the groups $H^2(U_2, \F_p)$ and  $H^2(U_1,\F_p)$ is trivial, then the claim is obvious. Suppose that both are non-trivial. So both are isomorphic to $\F_p$. By Corollary \ref{corestriction},  the corestriction map $\Cor \po H^2(U_1,\F_p) \lrar H^2(U_2, \F_p)$ is onto and so it is an isomorphism.   Since  the composition $\Cor \circ \Res$ is equal to the endomorphism that consists on multiplication by $|U_2: U_1|$, $\Cor \circ \Res=0$, and so $\Res$ should be the trivial map. \renewcommand\qedsymbol{$\diamond$} 
 \end{proof}
We are ready to prove in Claim \ref{cdsub} that $G$ is $p$-good at $\mathcal C$. Notice that the second conclusion of Claim \ref{cdsub} is inspired by Serre's profinite analogue \cite[Exercise 5(b)]{Ser97} of Strebel's theorem \cite{Str77} on infinite-index subgroups of $\PD_n$ groups.
   \begin{claim} \label{cdsub}
    Let $\mathcal C$ be pseudovariety of finite groups such that $\mathcal N_p\mathcal C=\mathcal C$.
Then any subgroup $H$ of $G$ of finite index is cohomologically $p$-good in $\mathcal C$.
In particular, $\cd_p(G_{\widehat {\mathcal C}})\leq 2$. 
Moreover, for any closed  subgroup $\mathbf H$ of $G_{\widehat {\mathcal C}}$ such that  $p^\infty $ divides the index $|G_{\widehat {\mathcal C}}:\mathbf H|$, it follows that $\cd_p(\mathbf H)\leq 1$.
 \end{claim}
  \begin{proof} Since $H$ is finitely generated and of cohomological dimension 2, we only have to check that the map  (\ref{pgoodC}) for  $H\lrar  H_{\widehat {\mathcal C}}$ is an isomorphism for dimension $n=2$,
 By Claim \ref{restriction}, the condition ($D_2$) of \cite[Section I.2.6]{Ser97} holds for  the embedding $H\hookrightarrow  H_{\widehat {\mathcal C}}$. Hence $H$  is cohomologically $p$-good in $\mathcal C$. For the second part, we observe that, if $\mathbf U$ is an open subgroup of $\mathbf H$,
 then $$\Hc ^k(\mathbf U, \F_p)=\varinjlim_{\mathbf U\leq \mathbf V\le_o \Gc} \Hc ^k(\mathbf V, \F_p)=\varinjlim_{\mathbf U\leq \mathbf V\le_o \Gc} H^k(\mathbf V\cap G, \F_p)=0,$$ where the second equality is a consequence of  the cohomological $p$-goodness that we have just proved.  Hence $\Hc ^k(\mathbf U, \F_p)=0$ if $k\geq 2$ by Claim \ref{restriction}. Thus, $\cd_p(\mathbf H)\leq 1$. \renewcommand\qedsymbol{$\diamond$}
 \end{proof}

 Claim \ref{cdsub} implies that if we lied under the assumptions of part (I) of Theorem \ref{Mel'2}, then its profinite completion $\hat G$ would satisfy $\cd_p(\hat G)=1$ for all primes $p$, and hence $\hat G$ would be projective \cite[Theorem 7.6.7]{Rib00}. This completes the proof of part (I). Before proving part (II) of Theorem \ref{Mel'2}, we first show that $G$ is 2-free. We proceed by contradiction. Let us suppose that $G$ is not 2-free.  Then, by \cite[Theorem 1.5 and Definition 6.5]{Lou22}, there exists a subgroup $K$ of $F$ of rank 2 containing $w$ such that the canonical homomorphism $P=K/\langle\!\langle w\rangle\!\rangle \lrar G$ is an embedding of a non-free group $P$ into $G$. From this assumption, we will now derive several claims that will lead to a contradiction.

 \begin{claim}\label{smallH_2}
 For any normal subgroup $H\leq G$ of finite index and every prime $p$, we have that   $\dim_{\F_p} H^2(H\cap P;\F_p)\leq 1$.
 \end{claim}  
 \begin{proof} 
     Let $\langle\!\langle w\rangle\!\rangle \leq U\lhd F$ be such that $H=U/\langle\!\langle w^F\rangle\!\rangle$. We take a transversal   $\{t_1, t_2,\ldots, t_{|K:U\cap K|}\}$ of $U\cap K$ in $K$, with $t_1=1$, and then we complete it to a   transversal $T=\{t_1,\ldots, t_{|F:U|}\}$ of $U$ in $F$.  We denoted by $n$ the rank of $F$. By the Nielsen--Schreier formula, $d(U)=1+(n-1)|F: U|$. By assumption, the group $H=U/\langle\!\langle w^{t_1},\ldots, w^{t_{|F:U|}}\rangle\!\rangle  $ is  one-relator. Let $m$ be such that $U/\langle\!\langle w^{t_1},\ldots, w^{t_{|F:U|}}\rangle\!\rangle  $ admits an $m$-generator one-relator presentation. By the multiplicativity of the Euler characteristic we have that $m=2+(n-2)|F: U|.$ From the fact that the $p$-abelianisation of $H$ is $m$-generated, it follows that the images of $w^{t_1},\ldots, w^{t_{|F:U|}}$ in $U/[U,U]U^p$ generate a subspace of dimension at least $d(U)-m=|F: U|-1$.  Henceforth, the images of $w^{t_1},\ldots, w^{t_{|K:U\cap K|}}$ in $U/[U,U]U^p$ generate  a subspace of dimension at least $|K:U\cap K|-1$. In particular,  the images of $w^{t_1},\ldots, w^{t_{|K:U\cap K|}}$ in $(U\cap K)/[U\cap K,U\cap K](U\cap K)^p$ generate a subspace of dimension at least $|K:U\cap K|-1$. Recall that the free group $K$ has rank two and so $d(U\cap K)=|K:U\cap K|+1$. 
We also observe that 
$H\cap P=U\cap (K\langle w^F\rangle)/\langle w^F\rangle=(U\cap K)\langle w^F\rangle/\langle w^F\rangle=U\cap K/\langle w^U\rangle$. 
From this, we obtain that    $\dim_{\F_p} H^1(H\cap P ;\F_p)\leq 2$. Again, by the multiplicativity of $\chi$,  $\chi(H\cap P )=|P: H\cap P|\, \chi(P)=0$. Thus $\dim_{\F_p} H^2(H\cap P;\F_p)\leq 1$. \renewcommand\qedsymbol{$\diamond$}
 \end{proof}

 \begin{claim}\label{b1n}
The group $P$ is isomorphic to $B(1, m)$. 
 \end{claim}
 \begin{proof}
 Let $\overline P$ be the closure of $P$ in $\widehat{G}$. Let $\mathbf W$ be an open subgroup of $\overline P$  and let $\mathbf U$ be an open normal subgroup of $\widehat{G}$ such that $\mathbf U\cap \overline P\leq \mathbf W$. We put $U=\mathbf U\cap G$ and $W=\mathbf W\cap P$.  Fix a prime $p$. 
 Observe that $U\cap P$ is a subgroup of $W$ and so,  by Corollary \ref{corestriction},   the corestriction map $H^2(U\cap P;\F_p)\lrar H^2(W;\F_p)$ is onto.
 Thus, we have that
 \begin{align*}
 \dim_{\F_p} \Hc ^1(\mathbf W;\F_p) & \leq \dim_{\F_p} H^1(W;\F_p)=1+\dim_{\F_p} H^2(W;\F_p)
 \\ &\leq 
 1+\dim_{\F_p} H^2(U\cap P;\F_p)\myle{Claim \ref{smallH_2}} 2.\end{align*}
 Hence the pro-$p$ Sylow subgroups of $\overline P$ are of bounded rank. We know  from Claim \ref{cdsub} that $\cd_p(\overline P)=1$ and hence that the pro-$p$ Sylows of $\overline{P}$ are free pro-$p$ groups \cite[Section 7.7]{Rib00}. Thus, in fact,  the 
 pro-$p$ Sylow subgroups of $\overline P$ are pro-$p$ cyclic. A profinite group all of whose pro-$p$ Sylows are pro-$p$ cyclic must be meta-cyclic (see \cite[Lemma 4.2.5]{Rib17}, this follows from the analogous result for finite groups due to Zassenhaus). 
 
 Since $P$ is embedded in $\overline P$, then $P$ is  meta-abelian.  
So we know that $P$ is non-cyclic, that it contains no non-abelian free subgroups and that it splits as a HNN extension $P=V*_{E, \theta}$ of a one-relator group $V$ along an injective map of a non-trivial free Magnus subgroup $\theta\po E\lrar V$. Henceforth, $E\cong \Z$. Moreover, if $E$ was a proper subgroup of $V$ and  $\theta$ was not surjective, then we could produce a non-abelian free subgroup of $P$ using Bass--Serre theory. It follows that $P\cong B(1, m).$  \renewcommand\qedsymbol{$\diamond$}   
 \end{proof}
 
 \begin{claim} \label{Uclaim}
 For any prime $p$, there exists  $\langle\!\langle K\rangle\!\rangle \leq U\lhd F$ such that $|F:U|=p$.
 \end{claim}
 \begin{proof}
 In the proof of Claim \ref{b1n} we have shown that the 
 pro-$p$ Sylow subgroups of $\overline P$ are pro-$p$ cyclic. Hence $p^{\infty}$ divides $|\widehat{G}/\overline {\langle\!\langle P\rangle\!\rangle} [\widehat{G}, \widehat{G}]|$. \renewcommand\qedsymbol{$\diamond$}
 \end{proof}
 In order to reach a contradiction, we consider two possible cases about $P$. Suppose first  that $P\cong B(1, m)$ with $m\ne 2$. Then, for any prime $p$ dividing $m-1$,  $w\in [K,K]K^p$.
  Consider the subgroup $U$ from Claim \ref{Uclaim} and name $H=U/\langle\!\langle w\rangle\!\rangle $, which has finite index in $G$. Since $w\in [U,U]U^p$, $$\dim_{\F_p} H^1(H;\F_p)=(n-1)p+1>(n-2)p+2.$$ But this contradicts Claim \ref{abfi}. Secondly, suppose that $P\cong B(1,2)$, there exists a subgroup $G_1$ of $G$ of finite index such that $P\cap  G_1\cong B(1, m)$ with $m\ne 2$. Applying a similar argument   as before we also reach a contradiction. We have proved the following. 
  \begin{claim} \label{2free} The group $G$ is 2-free.
  \end{claim}
 Now we move on to show part (II) of Theorem \ref{Mel'Intro}. We assume that there exists a subgroup $H$ of $G$ of finite index such that $H^2(H;\F_p)\cong \F_p$. Recall that a torsion-free group that is virtually isomorphic to a surface group must be a surface group itself by Kerckhoff's realisation theorem \cite{Ker83}. So we can assume that $H=G$. 
 
By \cite{Lin22}, it follows from Claim \ref{2free} that $G$ is hyperbolic and virtually special (in the sense of Haglund--Wise \cite{Hag08}). Again, by \cite{Hag08}, $G$ virtually embeds into a Right-angled Artin group, which are $\mathcal N_q$-RFRS by \cite{Kob20}. Henceforth, $G$ is virtually  $\mathcal N_q$-RFRS. Again, for notational convenience, we can assume that $G$ itself is $\mathcal N_q$-RFRS. Put $\mathcal C=\mathcal N_p(\mathcal A_f \mathcal N_q)$ and assume the notation of Section \ref{sect:rigidity}. We already have some of the assumptions that are required to apply  Theorem \ref{mainC0} to $G$ and hence conclude that $G$ is a surface group. We proved that $G$ is $p$-good at $\mathcal C$ in Claim \ref{cdsub}. It remains to show that $G_{\hat{\mathcal C}}$ is Poincar\'e Duality of dimension two at $p$. 

\begin{claim} \label{H2dem} We  have that $\Hc^2(\Gc, \Gcr)\cong \F_p$. 
\end{claim}
\begin{proof} By Proposition \ref{inverselimit}, we have an isomorphism
 \begin{equation*}
\begin{tikzcd}[cramped, sep=small]
     \Hc^2 (G_{\widehat{\mathcal C}};\Gcr)\ar[r]&
     \varprojlim_{G/N\in \mathcal C} \Hc^2 (G_{\widehat{\mathcal C}}; \F_p[G/N])
\end{tikzcd}
\end{equation*} 
Using Shapiro's lemma and the $p$-goodness of $G$ at $\mathcal C$, we can reformulate the above to get the isomorphism   \begin{equation} \label{H2Demeq}
\begin{tikzcd}[cramped, sep=small]
     \Hc^2 (G_{\widehat{\mathcal C}};\Gcr)\ar[r]&
     \varprojlim_{G/N\in \mathcal C} H^2 (N; \F_p),
\end{tikzcd}
\end{equation} 
where, given two open subgroups $N_1\leq N_2$ of $G$ in the pro-$\mathcal C$ topology,  the corresponding map $H^2 (N_1; \F_p)\lrar H^2(N_2; \F_p)$ that defines the inverse limit above is the corestriction map. Now we gather three observations. 
\begin{itemize}
    \item The corestriction map is always onto    (Corollary \ref{corestriction}). 
    \item Each group $H^2 (N; \F_p)$ is either trivial or $\F_p$ (because $N$ is one-relator).
    \item By assumption, $H^2(G; \F_p)\cong \F_p$.
\end{itemize}
From these, it is not hard to see that every  $H^2 (N; \F_p)$ must be isomorphic to $\F_p$, that the correstriction maps are isomorphisms and hence, by equation (\ref{H2Demeq}), that $\Hc^2 (G_{\widehat{\mathcal C}};\Gcr)\cong \F_p$.\renewcommand\qedsymbol{$\diamond$}
\end{proof}

The last ingredient will be to prove that $\Hc^1(\Gc, \Gcr)$ vanishes. This will be done in our last two claims. 
\begin{claim}\label{Demushkin}
For any subgroup $H\leq G$ of finite index, $H_{\widehat{\mathcal N_p}}$ is Demushkin. 
\end{claim}

\begin{proof} By the first part of Claim \ref{cdsub}, $H_{\widehat{\mathcal N_p}}$ is non-free one-relator pro-$p$ group. By Claim \ref{abfi}, any open subgroup $\UU\leq_o H_{\widehat{\mathcal N_p}}$ has the property that \[\dim_{\F_p} \Hc^1(\UU;\F_p)-2\leq (\dim_{\F_p} \Hc^1(H_{\widehat{\mathcal N_p}};\F_p)-2)|H_{\widehat{\mathcal N_p}}: \UU|.\] It follows from \cite{An73} (and also from  \cite{Dum83}) that $H_{\widehat{\mathcal N_p}}$ is Demushkin. \renewcommand\qedsymbol{$\diamond$}
\end{proof}
We do not really need to delve into the definition of Demushkin groups. However, for the purpose of Claim \ref{H1dem}, we recall that Demushkin groups are the Poincar\'e Duality pro-$p$ groups of dimension two (see \cite[Section I.4.5, Example 2]{Ser97}).
\begin{claim} \label{H1dem} We  have that $\Hc^1(\Gc, \Gcr)=0$. 
\end{claim}
\begin{proof}
    Let $\mathbf H$ be an open normal subgroup of $ \Gcr$. Then  $\mathbf H=\overline H$, where $H=\mathbf H\cap G$. Let  $\mathbf N_ H$ the kernel of the canonical map
$\mathbf H\lrar  H_{\widehat {\mathcal N_p}}$. Then we can write $G_{\hat{\mathcal C}}$ as the following inverse limit $$\Gc\cong \varprojlim_{\mathbf N\lhd_o \Gc} \Gc/\mathbf N_H.$$ This implies, by Proposition \ref{inverselimit}, that $$ \Hc^1(\Gc, \Gcr)\cong  
\varprojlim_{\mathbf N\lhd_o \Gc} \Hc^1(\Gc,\F_p\llbracket \Gc/\mathbf N_H \rrbracket ) .$$ 
Using the fact that $H$ has finite index in $G$ and  Shapiro's lemma we obtain  that $$\Hc^1(\Gc,\F_p\llbracket\Gc/\mathbf N_H \rrbracket )\cong \Hc^1(\mathbf H,\F_p\llbracket\mathbf H/\mathbf N_H \rrbracket )=\Hc^1(  H_{\widehat {\mathcal N_p}},\F_p\llbracket H_{\widehat {\mathcal N_p}} \rrbracket ).$$ 
Lastly, by Claim \ref{Demushkin}, $H_{\widehat {\mathcal N_p}} $ is Demushkin and hence the right-most cohomology group is trivial. It follows that $\Hc^1(\Gc, \Gcr)=0$. \renewcommand\qedsymbol{$\diamond$}
\end{proof}

By Claims \ref{cdsub}, \ref{H2dem} and \ref{H1dem}; we have that $\cd_p(\Gc)\leq 2$, $\Hc^2(\Gc, \Gcr)\cong \F_p$ and $\Hc^1(\Gc, \Gcr)=0$. Henceforth $\Gc$ is Poincar\'e Duality of dimension two over $\F_p$ at $p$ (recall Definition \ref{PD2profinite}). Finally, by Theorem \ref{mainC0}, $G$ is a surface group  and the proof of Theorem \ref{Mel'2} is complete.  \end{proof}

\bibliographystyle{amsalpha}

\bibliography{biblio}

\end{document}